\documentclass[12pt]{article}
\usepackage[T1]{fontenc}
\usepackage[latin9]{inputenc}
\usepackage{array}
\usepackage{verbatim}
\usepackage{float}
\usepackage{multirow}
\usepackage{amsthm}
\usepackage{amsmath}
\usepackage{amssymb}
\usepackage{graphicx}
\usepackage{setspace,color}
\usepackage{epstopdf}
\usepackage{esint}
\usepackage[authoryear]{natbib}
\usepackage{xargs}[2008/03/08]
\usepackage[textwidth=6.5in,textheight=9in]{geometry}
\usepackage[subrefformat=parens]{subfig}
\usepackage{changepage}

\makeatletter

\providecommand{\tabularnewline}{\\}

\newtheorem{theorem}{Theorem}

\newenvironment{lyxlist}[1]
{\begin{list}{}
{\settowidth{\labelwidth}{#1}
 \setlength{\leftmargin}{\labelwidth}
 \addtolength{\leftmargin}{\labelsep}
 }}
{\end{list}}

\@ifundefined{date}{}{\date{}}
\usepackage{ifpdf} 
\ifpdf 

 \IfFileExists{lmodern.sty}{\usepackage{lmodern}}{}

\fi 

\usepackage{mathrsfs}
\usepackage{enumitem}
\usepackage{cancel}

\setenumerate{label= \textup{(}\textit{\roman*}\textup{)},leftmargin=1.6pc,topsep=5pt}

\usepackage{ifthen}
\def\boxit#1{\vbox{\hrule\hbox{\vrule\kern6pt
          \vbox{\kern6pt#1\kern6pt}\kern6pt\vrule}\hrule}}

\def\LyX{\texorpdfstring{%
  L\kern-.1667em\lower.25em\hbox{Y}\kern-.125emX\@}
  {LyX}}

\makeatother

\usepackage[left,displaymath,mathlines,pagewise]{lineno}
\usepackage{xcolor}

\renewcommand\linenumberfont{\normalfont\scriptsize\sffamily\color{red}}
\global\long\def\real{\mathbb{R}}

\global\long\def\myvec#1{\mathbf{#1}}

\global\long\def\probop#1{\mathbf{#1}}

\global\long\def\expect{\probop E}
\global\long\def\real{\mathbb{R}}

 \global\long\def\diff{\operatorname{d}}
 \global\long\def\ssetvar#1{\mathcal{#1}}
 
 \newcommandx\finiteseq[5][usedefault, addprefix=\global, 1={,}, 2=1, 3=2, 4=n]{#5_{#2}#1#5_{#3}#1\ldots#1#5_{#4}}

\newcommandx\infseq[4][usedefault, addprefix=\global, 1={,}, 2=1, 3=2]{#4_{#2}#1#4_{#3}#1\ldots}
 \newcommandx\limitsup[2][usedefault, addprefix=\global, 1=n, 2=\infty]{\underset{#1\rightarrow#2}{\lim\sup}}
 \newcommandx\limitinf[2][usedefault, addprefix=\global, 1=n, 2=\infty]{\underset{#1\rightarrow#2}{\lim\inf}}

 \newcommandx\probspace[3][usedefault, addprefix=\global, 1=\Omega, 2=\ssetvar E, 3=P]{(#1,#2,#3)}
 \newcommandx\topospace[2][usedefault, addprefix=\global, 1=X, 2=\tau]{(#1,#2)}
 \newcommandx\measurespace[3][usedefault, addprefix=\global, 1=X, 2=\ssetvar M, 3=\mu]{(#1,#2,#3)}
 \newcommandx\measurablespace[2][usedefault, addprefix=\global, 1=X, 2=\ssetvar M]{(#1,#2)}

{}

 \global\long\def\define{\overset{\text{def}}{=}}
 
 \newcommandx\basis[1][usedefault, addprefix=\global, 1=B]{\ssetvar{#1}}

{} \global\long\def\converge{\rightarrow}

 \newcommandx\conexp[2][usedefault, addprefix=\global, 1=\xi, 2=\ssetvar G]{\expect(#1|#2)}

\global\long\def\bdU{\myvec U}
\global\long\def\bdy{\myvec y}

\global\long\def\bdb{\myvec b}

\global\long\def\bdV{\myvec V}
\global\long\def\bdW{\myvec W}

\global\long\def\diffop{\mathcal{D}}

\global\long\def\penalized{\mathrm{smooth}}

\global\long\def\bdB{\myvec B}
\global\long\def\trans{T}

\title{A Smooth and Locally Sparse Estimator for Functional Linear Regression
via Functional SCAD Penalty}
\author{Zhenhua Lin$^{1}$,
Jiguo Cao$^{2}$, Liangliang Wang$^{2}$ and Haonan Wang$^{3}$}

\begin{document}

\maketitle
\noindent
{\small $^{1}$Department of Statistical Sciences, University of Toronto, Toronto, ON, Canada} \\
{\small $^{2}$Department of Statistics and Actuarial Science, Simon Fraser University,   Burnaby, BC, Canada}\\
{\small $^{3}$Department of Statistics, Colorado State University, Fort Collins, CO, U.S.A.}

\begin{abstract}
In this paper, we propose a new regularization technique called ``functional SCAD''. We then combine this technique with the smoothing spline method to develop a
smooth and locally sparse (i.e., zero on some sub-regions) estimator
for the coefficient function in functional linear regression. The functional SCAD has a nice shrinkage property that enables
our estimating procedure to identify the null subregions of the coefficient
function without over shrinking the non-zero values of the coefficient function. Additionally, the smoothness
of our estimated coefficient function is regularized by a roughness penalty rather than by controlling the number of knots. Our method is more theoretically sound and is computationally
simpler than the other available methods. An asymptotic analysis shows that our estimator is consistent and can identify the null region with the probability tending to one. Furthermore, 
simulation studies show that our estimator has superior numerical
performance. Finally, the practical merit of our method is demonstrated on two real applications.
\end{abstract}

\textit{keywords: }
B-spline basis, null region, smoothly clipped
absolute deviation.%

\section{Introduction\label{sec:introduction}}

In this paper, we study the problem of estimating the coefficient function $\beta(t)$ in the following functional linear regression (FLR) model
\begin{equation}
Y_{i}=\mu+\int_{0}^TX_{i}(t)\beta(t)\diff t+\varepsilon_{i},  \label{eq:flr-model}
\end{equation}
where each $Y_{i}$ is a scalar response, $\mu$ is the grand mean,
$X_{1}(t), \ldots, X_n(t)$ are independent realizations of an unknown
centered stochastic process $X(t)$ defined on the domain $[0,T]$,
and $\varepsilon_{i}\overset{i.i.d.}{\sim}N(0,\sigma_{\varepsilon}^{2})$, $i = 1, \ldots, n$.
The coefficient function $\beta(t)$ is also called the regression weight
function: $\beta(t)$ weights the cumulative contribution of a functional predictor
$X_{i}(t)$ to the corresponding response $Y_{i}$.
It is of practical interest to know where $X_i(t)$ has no effect on $Y_i$,
where $X_i(t)$ has a significant effect on $Y_i$, and the effect  magnitude. Statistically, this paper focuses on 
 identifying the null subregions
of $\beta(t)$ (i.e., the subregions on which $\beta(t)=0$), and simultaneously estimating $\beta(t)$ on its non-null
subregions.

Historically, functional linear regression originates from ordinary
linear regression with a large number of predictors. When the number
of predictors is very large, estimators produced by ordinary least
squares exhibit excessive variability, and hence perform poorly 
on prediction. Researchers have developed many approaches to rectify
the issue, such as partial least squares (PLS), principal components
regression (PCR), and ridge regression, etc. Unfortunately, when the predictors were
actually discrete observations of some continuous process, \citet{Hastie1993}
pointed out that none of these methods directly made use of the spatial
information, i.e., the correlation and the order between predictors.
They suggested that it would be more natural to regularize the variability
by directly restricting the coefficient vector to be smooth.
They also prototyped the ideas of penalized least squares and smoothing splines. These ideas were then explored more thoroughly in
\citet{Ramsay1997}.

Following \citet{Hastie1993} and \citet{Ramsay1997}, more theoretical and practical treatments of functional linear regression emerged. For example, \citet{Cardot2003}
studied some theoretical aspects of FLR, such as identifiability,
existence and unicity of the estimator of $\beta(t)$. They also proposed
an estimator based on smoothing splines with a roughness
penalty. The method developed by \citet{Crambes2009} is also based
on smoothing splines, but it includes an unusual term in
the roughness penalty to ensure the existence of their estimator without
additional assumptions on sample curves. Both works used 
B-spline bases  and provided
convergence rates of the resulting estimators  in terms of a norm 
 induced by
the covariance operator of the process $X(t)$. Instead of B-spline approximation,
\citet{Li2007} used a penalized Fourier basis to estimate $\beta(t)$, and derived
an $L^{2}$-norm convergence rate for their estimator. While
 \citet{Hastie1993}  motivated  the aforementioned works  to use penalized least squares, some researchers also used spectral decomposition of the covariance function of $X(t)$ to estimate the coefficient function. For instance, \citet{Cardot2003} 
 proposed to obtain an initial estimate through principal components
regression and then smooth it by means of B-spline approximation.
In \citet{Cai2006} and \citet{Hall2007}, $\beta(t)$ was estimated
based on the spectral decomposition as well.

Besides the classic task of estimating $\beta(t)$, researchers have also
studied some other aspects of FLR. For example, \citet{Cardot2007}
investigated the convergence rates of some estimators in the weak topology
of $L^{2}([0,T])$. \citet{Preda2007} and \citet{Yuan2010} explored
FLR from the point of view of reproducing kernel Hilbert space. \citet{Muller2005}
studied generalized functional linear models. \citet{Yao2005} extended
the scope of FLR to the sparse longitudinal data. \citet{Fan2000}
and \citet{Lian2012} studied FLR with a functional response.

Recently,
\citet{James2009} recognized the importance of the interpretability
of estimated $\hat{\beta}(t)$ and proposed a so-called ``FLiRTI'' approach
that intended to produce a more interpretable estimator. In order to obtain  $\hat{\beta}(t)$ which is zero on  subregions, FLiRTI penalizes the function values and  derivatives  at some grid points using the L1 regularization method. This method can successfully shrink small values of $\hat{\beta}(t)$. Unfortunately, FLiRTI cannot warrant that the estimated $\hat{\beta}(t)$ is zero on some subregions for the following two reasons: 1) in
order to have  $\hat{\beta}(t)$ to be zero on some subregions, we need to force the function
value at the group of grid points in those subregions to be zero \textit{simultaneously}; 2) even 
we force  a function and its first several derivatives  to be zero at a point by L1 regularization, 
it may still nonzero  at a subregion around the point: a smooth function is zero in a subregion around a point if and only if its value and \textit{all} of its derivatives are zero at that point.  
 As pointed out by \citet{Zhou2013}, the FLiRTI estimate tends to have a large variation: 
a small grid size  causes unstable numerical solution of the
LASSO method;  but a large set of dense grid points makes FLiRTI tend to overparametrize the model.

\citet{Zhou2013} also proposed a method for simultaneously identifying
the null subregions of $\beta(t)$ and estimating $\beta(t)$ on
its non-null subregions. The method has two stages: in the first
stage, a Dantzig selector is used to obtain an initial estimate of the
null subregions, and in the second stage, the estimated null subregions
are refined. The estimate of $\beta(t)$
on the non-null subregions is produced via a group SCAD penalty proposed
by \citet{Wang2007}. 
Selection of several tuning parameters in each stage not only  increases the estimation variability and computational complexity but also makes the method
difficult to implement, which  limits the applications of this method.

In this paper, we first propose a new regularization technique called
``functional SCAD'' (fSCAD for short), and then combine it with smoothing splines to develop a one-stage procedure, called \textit{SLoS} (\textit{S}mooth
and \textit{Lo}cally \textit{S}parse). This procedure simultaneously identifies
the null subregions of $\beta(t)$ and produces a smooth estimate
of $\beta(t)$ in the non-null subregions. Here,  ``locally sparse''  means a curve is zero on some subregions of its domain \citep{Tu2012}. The fSCAD
can be viewed as a functional generalization of the ordinary SCAD
proposed in \citet{Fan2001}. Its nice shrinkage property 
enables the SLoS estimating procedure to identify the null subregions of $\beta(t)$
without over shrinking the non-zero values of $\beta(t)$ on the non-null subregions. In addition, we employ a roughness
penalty to regularize the smoothness of the SLoS estimator of $\beta(t)$. Compared 
to existing methods in the literature, our SLoS method has two distinct
features. First, unlike the FLiRTI pointwise penalization method, our fSCAD regularizes
the ``overall magnitude'' of the estimated coefficient function $\hat{\beta}(t)$. Second, unlike the two-stage procedure
in \citet{Zhou2013}, our SLoS method combines  fSCAD and smoothing splines together in a single optimization objective function. By solving this optimization problem, we are able to 
produce an estimate of the null subregions, as well as a smooth estimate
of $\beta(t)$ in its non-null subregions simultaneously in one single
step. These two features make our estimating procedure theoretically sounder and computationally simpler.

The rest of the paper is organized as follows. In Section \ref{sec:method}
we propose the fSCAD regularization technique and the SLoS estimator. In Section \ref{sub:fscad}, we present the functional SCAD penalty and discuss its shrinkage property. We then provide a brief background of smoothing spline method in Section \ref{sub:smoothing-spline}. In Section \ref{sub:Locally-Sparse-Smooth}, we propose the SLoS estimator for the coefficient function $\beta(t)$. In Section \ref{sub:Extension}, we extend our method to the case where there are multiple functional predictors.
The asymptotic properties of the SLoS estimator are provided in Section \ref{sec:asymptotic}.
In Section \ref{sec:Simulations} we evaluate the numerical performance
of our estimator via simulation studies. In Section \ref{sec:Applications}
we apply the SLoS method to study two real datasets. The paper
is concluded in Section \ref{sec:Conclusions}. Proofs of theorems
are collected in a supplementary file. The Matlab codes for running the application and simulation studies can be downloaded at the website \textsf{http://people.stat.sfu.ca/$\sim$\,cao/Research/FLR/}.

\section{Methodology\label{sec:method}}

\subsection{Functional SCAD}\label{sub:fscad}

To motivate our functional SCAD, we first briefly review the ordinary SCAD penalty proposed by \citet{Fan2001}
in the setting of regression with multiple scalar predictors. Here,
we use ordinary linear regression as an example to illustrate
the SCAD technique. An ordinary linear model with one response variable and $J$ covariates can be expressed as
\begin{eqnarray}\label{eqn:lm}
y_{i}=\mu+\sum_{j=1}^{J}b_{j}z_{ij}+\varepsilon_{i},
\end{eqnarray}
where $\varepsilon_{i}\overset{i.i.d.}{\sim}N(0,\sigma^{2}), i = 1,\ldots, n$. Using
the least-squares function as a loss function, the SCAD estimator of
$b_{1},b_{2},\ldots,b_{J}$ is defined as
\begin{equation}
\hat{\myvec b}^{\mathrm{scad}}=\underset{\myvec b\in\real^{J}}{\arg\min}\left\{ \frac{1}{n}\sum_{i=1}^{n}\left(y_{i}-\mu-\sum_{j=1}^{J}b_{j}z_{ij}\right)^2+\sum_{j=1}^{J}p_{\lambda}(|b_{j}|)\right\} ,\label{eq:olr-scad}
\end{equation}
where the penalty function $p_{\lambda}(\cdot)$ is the SCAD penalty
function of \citet{Fan2001}. It is defined on $[0,+\infty]$ as
\[
p_{\lambda}(u)=\left\{ \begin{array}{ll}
\lambda u & \textrm{if}\,\,\,0\leq u\leq\lambda\\
-\frac{u^{2}-2a\lambda u+\lambda^{2}}{2(a-1)} & \textrm{if}\,\,\,\lambda<u<a\lambda\\
\frac{(a+1)\lambda^{2}}{2} & \textrm{if}\,\,\, u\geq a\lambda,
\end{array}\right.
\]
where a suggested value for $a$ is $3.7$ according to \citet{Fan2001},
and $\lambda$ is a tuning parameter varying with sample size. Additionally, 
when the number of covariates is fixed, the SCAD penalty enjoys the so-called
oracle property: it is able to identify the true sub-model with probability
tending to one, and meanwhile produce an asymptotically normal estimate
for each non-zero variable (\citet{Fan2001}). In other words, the estimator $\hat{\myvec b}^{\mathrm{scad}}$
performs as well as if the true sub-model was known, i.e., we know which coefficients $b_{j}$'s are zero in advance. Even when the number of covariates grows as the sample size increases, $\hat{\myvec b}^{\mathrm{scad}}$
still enjoys the oracle property under certain regularity conditions, as discussed in \citet{Fan2004} and \citet{Fan2011}.

The ordinary linear model (\ref{eqn:lm}) changes to the functional linear model (\ref{eq:flr-model}) if the summation $\sum_{j=1}^{J}b_{j}z_{ij}$ is replaced by the integral $\int_{0}^TX_{i}(t)\beta(t)\diff t$. Similarly, for the SCAD penalty term $\sum_{j=1}^{J}p_{\lambda}(|b_{j}|)$ in
(\ref{eq:olr-scad}), if we normalize it by $\frac{1}{J}$, then  in the functional space, it is generalized to an integral $\frac{1}{T}\int_{0}^Tp_{\lambda}(|\beta(t)|)\diff t$.
In this paper, we call 
\begin{eqnarray}\label{fSCAD}
\mathcal{L}(\beta)\define\frac{1}{T}\int_{0}^Tp_{\lambda}(|\beta(t)|)\diff t
\end{eqnarray}
 the fSCAD  penalty.

We have introduced the fSCAD as an analogy to the ordinary SCAD. From
that point of view, the fSCAD might be viewed as a functional generalization
of the ordinary SCAD. Now, to gain insight into its shrinkage nature,
we shall present a different way of introducing it in the setting
of locally sparse modeling. For example, suppose $\beta(t)$ is locally
sparse and we want to have a locally sparse estimate of $\beta(t)$.
A natural idea is to divide the domain into many small subintervals
and then penalize the overall magnitude of $\beta(t)$ on each
subinterval to shrink the estimated $\hat{\beta}(t)$ towards zero on those subintervals where the true $\beta(t)$ is zero. A possible measure of overall magnitude
could be a normalized $L^{2}$ norm of $\beta(t)$ on a subinterval,
and a possible penalty function could be the ordinary SCAD penalty.
The following remarkable result shows that this idea of locally sparse
modeling actually leads to our fSCAD regularization.
\begin{theorem}
\label{thm:penalty-term}Let $0=t_{0}<t_{1}<\cdots<t_{M}=T$ be an
equally spaced sequence in the domain $[0,T]$, and $\beta_{[j]}(t)$
denote the restriction of $\beta(t)$ on the subinterval $[t_{j-1},t_{j}]$,
i.e., $\beta_{[j]}(t)=\beta(t)$ for $t\in[t_{j-1},t_{j}]$ and zero
elsewhere. If $\beta(t)$ is continuous, then for any $q\geq 1$,
\begin{equation}
\frac{1}{T}\int_{0}^Tp_{\lambda}(|\beta(t)|)\diff t=\lim_{M\converge\infty}\frac{1}{M}\sum_{j=1}^{M}p_{\lambda}\left(M^{\frac{1}{q}}T^{-\frac{1}{q}}\|\beta_{[j]}\|_{q}\right),\label{eq:fscad-alt-def}
\end{equation}
where $\|\beta_{[j]}\|_{q}\define\left(\int_{0}^T|\beta_{[j]}(t)|^{q}\diff t\right)^{1/q}$.
\end{theorem}
In Theorem \ref{thm:penalty-term}, the normalized $L^q$ norm $M^{\frac{1}{q}}T^{-\frac{1}{q}}\|\beta_{[j]}\|_{q}$ plays the role of measuring ``overall magnitude'' of $\beta(t)$ over the subinterval $[t_{j-1},t_j]$. The identity (\ref{eq:fscad-alt-def}) shows that, the fSCAD of $\beta(t)$
is indeed the limit of the average SCAD penalty on the ``overall magnitude''
of $\beta(t)$ over each small subinterval $[t_{j-1},t_{j}]$. This
connection to the ordinary SCAD regularization sheds light on the
shrinkage nature of fSCAD. Roughly speaking, with fSCAD regularization,
Theorem \ref{thm:penalty-term} implies the following: 1) when $\beta(t)$ is overall very small on
$[t_{j-1},t_{j}]$ (i.e., $M^{\frac{1}{q}}T^{-\frac{1}{q}}\|\beta_{[j]}\|_{q}$
is very small), then the SCAD penalty $p_{\lambda}\left(M^{\frac{1}{q}}T^{-\frac{1}{q}}\|\beta_{[j]}\|_{q}\right)$
will shrink $\beta(t)$ towards zero identically over $[t_{j-1},t_{j}]$. This 
 indicates that fSCAD is able to choose a locally sparse estimate;
and 2) when $\beta(t)$ has significant overall  magnitude on $[t_{j-1},t_{j}]$
(i.e., $M^{\frac{1}{q}}T^{-\frac{1}{q}}\|\beta_{[j]}\|_{q}$ is big enough),
then $\beta(t)$ does not get penalized on the subinterval. This 
means fSCAD can avoid over shrinking $\beta(t)$.

Theorem \ref{thm:penalty-term} also provides a practical way to approximately
compute the fSCAD penalty: we choose a large $M$ and then
approximate $\frac{1}{T}\int_{0}^Tp_{\lambda}(|\beta(t)|)\diff t$ by
$\frac{1}{M}\sum_{j=1}^{M}p_{\lambda}\left(M^{\frac{1}{q}}T^{-\frac{1}{q}}\|\beta_{[j]}\|_{q}\right)$.
When B-spline expansion (introduced in the next section) is used to
approximate $\beta(t)$, we might take $q=2$ because it is relatively
easy to compute the $L^{2}$ norm of a B-spline function. Moreover,
for a B-spline basis, each $\|\beta_{[j]}\|_{2}$ only involves a few
basis coefficients due to the compact support property of B-spline
basis functions. In this particular setting, the approximated fSCAD can be viewed
as a generalized SCAD with a diverging number of parameters.

The shrinkage feature of
fSCAD brings it numerous potential applications in locally sparse
modeling. For example, in Section \ref{sub:Locally-Sparse-Smooth},
we use it to derive our SLoS estimator of the coefficient function
in functional linear regression. Now we begin to discuss the shrinkage nature of the fSCAD in a more
precise way. In other words, we will argue that the fSCAD shrinks the values of the estimator $\hat{\beta}(t)$ towards
zero on the null subregions of $\beta(t)$, but does not over shrink
the estimate on the non-null subregions.  First of all, since
\[
\left(\min_{t\in[t_{j-1},t_{j}]}|\beta_{[j]}(t)|^{q}\right)\frac{T}{M}\leq\int_{t_{j-1}}^{t_{j}}|\beta_{[j]}(t)|^{q}\diff t=\|\beta_{[j]}\|_{q}^{q}\leq\|\beta_{[j]}\|_{\infty}^{q}\frac{T}{M},
\]
where $\|\cdot\|_\infty$ denotes the supremum norm. Immediately, we have
\begin{equation}
\min_{t\in[t_{j-1},t_{j}]}|\beta_{[j]}(t)|\leq M^{\frac{1}{q}}T^{-\frac{1}{q}}\|\beta_{[j]}\|_{q}\leq\|\beta_{[j]}\|_{\infty}\leq\|\beta\|_{\infty}.\label{eq:penalized-term-bound}
\end{equation}
Let $H(\beta_{[j]})$ denote $M^{\frac{1}{q}}T^{-\frac{1}{q}}\|\beta_{[j]}\|_{q}$
and suppose $\hat{\beta}(t)$ is a consistent continuous estimator
of $\beta(t)$. For any fixed $t\in [0,T]$, if $\beta(t)\neq0$, then
by the continuity of $\beta(t)$ we have $|\beta(t)|>\epsilon$ on
a small neighborhood $\mathcal{N}_{\epsilon}(t)$ of $t$ for some
$\epsilon>0$. When $M$ is sufficiently large, the subinterval $[t_{j-1},t_{j}]$
containing $t$ is inside $\mathcal{N}_{\epsilon}(t)$. Then with
probability tending to one, we have $|\hat{\beta}_{j}(s)|\geq\epsilon/2$
for all $s\in\mathcal{N}_{\epsilon}(t)$ and hence $H(\hat{\beta}_{j})\geq\epsilon/2$.
As $\lambda_{n}\converge0$, $H(\hat{\beta}_{j})>a\lambda_{n}$ with
probability tending to one. This indicates that the consistent estimator
$\hat{\beta}(t)$ does not get penalized for its values at $t$. On
the other hand, if $\beta(t)=0$ on a small interval $\mathcal{N}_{0}$,
then $|\hat{\beta}(t)|\converge0$ on $\mathcal{N}_{0}$ with probability
tending to one. When $M$ is sufficiently large, one or more subintervals
$[t_{j-1},t_{j}]$ are inside $\mathcal{N}_{0}$. Then $H(\hat{\beta}_{j})\converge0$
in probability. By choosing an appropriate $\lambda_{n}$, the penalty
$p_{\lambda_{n}}(H(\hat{\beta}_{j}))$ grows with $H(\hat{\beta}_{j})$
in the rate $\lambda_{n}$, and hence forces $\hat{\beta}(t)$ to
become identically zero on $[t_{j-1},t_j]$.

\subsection{Smoothing Spline Method}\label{sub:smoothing-spline}

\begin{figure}[htbp]
\begin{centering}
\includegraphics[scale=0.7]{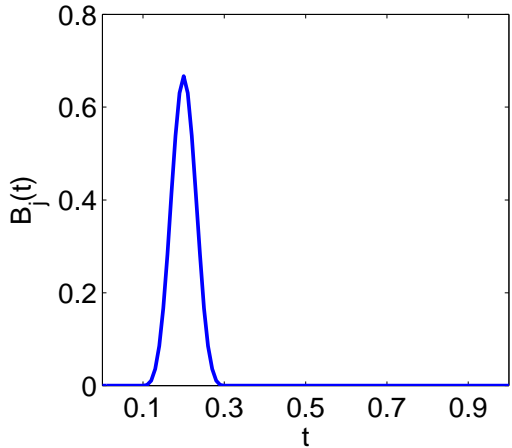}
\par\end{centering}
\caption{One example of a cubic B-spline basis function, which is only nonzero in [0.1,0.3]. The corresponding cubic B-spline basis function system is defined with 21 equally spaced knots in [0,1].}
\label{fig:example-b-spline}
\end{figure}

To prepare for the introduction of our SLoS estimator, in this section
we briefly review the smoothing spline method for estimating
$\beta(t)$. Recall that B-spline basis functions are defined by their order, as well as 
a sequence of knots. Suppose we set their order to $d+1$ and place $M+1$
equally spaced knots $0=t_{0}<t_{1}<\cdots<t_{M}=T$ in the domain
$[0,T]$ to define the B-spline basis functions. Over each subinterval, each B-spline
basis function $B_{j}(t)$ is a polynomial of the degree $d$. Moreover,
each B-spline basis function is nonzero over no more than $d+1$ consecutive
subintervals. When $M$ is large, each B-spline basis function is
only nonzero on a very short subregion. This property is called the
compact support property, which is critical for efficient computation
and makes B-spline bases very popular in functional data analysis.
This property is also very important for us to practically compute
the SLoS estimator. For example, Figure \ref{fig:example-b-spline} shows one of 23 cubic B-spline basis functions defined with 21 equally spaced knots. This basis
function is only nonzero on the interval {[}0.1,0.3{]}. Please refer to \citet{deBoor2001} for more details about B-spline basis functions.

Let $\mathcal{S}_{dM}$ denote the linear space spanned by the B-spline
basis functions $\{B_{j}(t):j=1,2,\ldots,M+d\}$  defined above.  
A modern method of estimating the coefficient function $\beta(t)$ 
in the functional linear model (\ref{eq:flr-model}) is to choose
$\beta(t) \in \mathcal{S}_{dM}$ which minimizes the sum of
squared residual errors. However, this method usually produces an estimator
that exhibits excessive variability, especially when $M$ is relatively
large. A popular approach to obtain a smooth estimator for $\beta(t)$ is to add regularization
through roughness penalty. For example, a smooth estimator of $\beta(t)$,
proposed in \citet{Cardot2003}, is defined as
\[
\hat{\beta}^{\penalized}=\underset{\beta\in\mathcal{S}_{dM}}{\arg\min}\left\{ \frac{1}{n}\sum_{i=1}^{n}\left[Y_{i} - \mu - \int_{0}^TX_{i}(t)\beta(t)\diff t\right]^{2}+\gamma\|\diffop^{m}\beta\|^{2}\right\},
\]
where $\diffop^{m}$ is the $m$th order
differential operator with $m\leq d$, and $\|\cdot\|$ denotes the $L^2$ norm of a function.  The parameter $\gamma\geq0$
is a tuning parameter varying with the sample size $n$. Note that
by setting $\gamma=0$, the above smooth estimator of $\beta(t)$ is reduced to the
classical ordinary least-squares estimator. This estimator is very rough when $M$ is large, while as $\gamma\converge\infty$,
the estimator $\hat{\beta}^{\penalized}$ converges to a linear
function. Thus, $\gamma$ serves as a smoothing parameter, which controls
the degree of smoothing of $\hat{\beta}^{\penalized}$.

Although regularization via roughness penalty is able to produce a
smooth estimator, it does not yield a locally sparse estimator.
In other words, it is not able to identify the null subregions of
$\beta(t)$. In the next section, we add the fSCAD to the fitting criterion to produce a smooth and meanwhile locally sparse estimate of $\beta(t)$.

\subsection{The SLoS Estimator\label{sub:Locally-Sparse-Smooth}}

We now combine fSCAD and smoothing splines together to obtain a smooth and locally
sparse estimator for $\beta(t)$ in the functional linear model (\ref{eq:flr-model}), which is called the SLoS estimator in this paper. 

The coefficient function $\beta(t)$ along with the constant intercept $\mu$ is estimated by minimizing $Q(\beta, \mu)$ defined as follows:
\begin{equation}
 Q(\beta, \mu) = \frac{1}{n}\sum_{i=1}^{n}\left[Y_{i}-\mu-\int_{0}^TX_{i}(t)\beta(t)\diff t\right]^{2}+\gamma\|\diffop^{m}\beta\|^{2}+\frac{M}{T}\int_{0}^Tp_{\lambda}(|\beta(t)|)\diff t,\label{eq:slos}
\end{equation}
where the minimization is taken over all $\beta(t)\in\mathcal{S}_{dM}$ and $\mu\in\real$. In the fitting criterion, $Q(\beta, \mu)$, we combine two penalty terms in a single optimization criterion: the
roughness penalty $\gamma\|\diffop^{m}\beta\|$ regularizes the
smoothness of $\beta(t)$, and the fSCAD penalty term $\frac{M}{T}\int_{0}^Tp_{\lambda}(|\beta(t)|)\diff t$ regularizes the local sparseness of $\beta(t)$. 

In the fitting criterion \eqref{eq:slos}, we scale the fSCAD penalty up by the factor $M$, where $M$ is the number of subintervals segmented by the knots of the B-spline basis functions for $\beta(t)$. This is necessary for the estimated coefficient function $\hat{\beta}(t)$ to enjoy the theoretical properties presented in Section 3. Also, according to Theorem 1,  the fSCAD penalty
term can be approximated by
\begin{eqnarray}\label{eqn:fSCADapprox}
\frac{M}{T}\int_{0}^Tp_{\lambda}(|\beta(t)|)\diff t\approx\sum_{j=1}^{M}p_{\lambda}\left(\frac{\|\beta_{[j]}\|_{2}}{\sqrt{T/M}}\right) 
= \sum_{j=1}^{M}p_{\lambda}\left(\sqrt{\frac{M}{T}\int_{t_{j-1}}^{t_{j}}\beta^{2}(t)\diff t}\,\right).
\end{eqnarray}
Thus, the fSCAD penalty term $\frac{M}{T}\int_{0}^Tp_{\lambda}(|\beta(t)|)\diff t$  actually represents the penalty on the overall magnitude of $\beta(t)$ of all $M$ subintervals. Based on this observation, we emphasize that the fSCAD penalty regularizes $\beta(t)$ in the way that it penalizes the overall magnitude of $\beta(t)$ on each subinterval, rather than the pointwise penalization fashion adopted in FLiRTI \citep{James2009}. This ensures that if the fSCAD regularization forces the overall magnitude of $\beta(t)$ on a subinterval to be zero, then the estimated $\hat{\beta}(t)$ is identically zero almost everywhere on that subinterval. Therefore, our method overcomes the shortcoming of the FLiRTI \citep{James2009}. Moreover, our  one-stage fitting procedure  is  computationally simpler than the two-stage method proposed by \citet{Zhou2013}.

Below we shall show how to compute the SLoS estimator $\hat{\beta}(t)$
in practice. To simplify the notation and manifest our main idea,
we assume $\mu=0$, and redefine the function $Q(\beta) = Q(\beta, \mu=0)$. The case that $\mu\neq0$ is considered at the end
of this section. We estimate $\beta(t)$ as a linear combination of B-spline basis functions:
\begin{eqnarray*}
\beta(t) = \sum_{k=1}^{M+d} b_k B_k(t) = \bdB^{\trans}(t)\bdb
\end{eqnarray*}
where $\bdB(t) = (B_{1}(t),B_{2}(t),\ldots,B_{M+d}(t))^{\trans}$ is the vector of B-spline basis functions, and $\bdb=(b_{1},b_{2},\ldots,b_{M+d})^{\trans}$
is the corresponding vector of coefficients. Let $\bdU$
be an $n\times(M+d)$ matrix with entries $u_{ij}=\int_{0}^TX_{i}(t)B_{j}(t)\diff t$
for $i=1,2,\ldots,n$, $j=1,2,\ldots,(M+d)$. Let $\bdy=(Y_{1},Y_{2},\ldots,Y_{n})^{T}$.
Then the first term of $Q(\beta)$ in (\ref{eq:slos}) is expressed in the matrix notation
as
\begin{equation}
\frac{1}{n}\sum_{i=1}^{n}\left[Y_{i}-\int_{0}^TX_{i}(t)\beta(t)\diff t\right]^{2}=\frac{1}{n}(\bdy-\bdU\bdb)^{T}(\bdy-\bdU\bdb).\label{eq:ls-part}
\end{equation}
Let $\bdV$ be an $(M+d)\times(M+d)$ matrix with entries $v_{ij}= \int_0^T \big(\frac{d^mB_{i}(t)}{dt^m}\frac{d^mB_{j}(t)}{dt^m}\big) dt $
for $1\leq i,j\leq M+d$. Then the second term of $Q(\beta)$ in
(\ref{eq:slos}) is
\begin{equation}
\gamma\|\diffop^{m}\beta\|^{2}=\gamma\bdb^{T}\bdV\bdb.\label{eq:roughness-part}
\end{equation}

Next, we turn to the third term of $Q(\beta)$ in \eqref{eq:slos}. According to Theorem \ref{thm:penalty-term}, we approximate the fSCAD penalty
term $\frac{M}{T}\int_{0}^Tp_{\lambda}(|\beta(t)|)\diff t$ as (\ref{eqn:fSCADapprox}),
and
\[
\int_{t_{j-1}}^{t_{j}}\beta^{2}(t)\diff t=\bdb^{T}\bdW_{j}\bdb\,,
\]
where $\bdW{}_{j}$ is an $(M+d)$-by-$(M+d)$ matrix with entries
$w_{uv}=\int_{t_{j-1}}^{t_{j}}B_{u}(t)B_{v}(t)\diff t$ if $j\leq u,v\leq j+d$
and zero otherwise. The SCAD penalty function $p_{\lambda}(\cdot)$
might be approximated by the local quadratic approximation (LQA) proposed
in \citet{Fan2001} as follows, since it makes the computation simpler. When $u\approx u^{(0)}$, the LQA
of SCAD function $p_{\lambda}(u)$ is
\[
p_{\lambda}(|u|)\approx p_{\lambda}(|u^{(0)}|)+\frac{1}{2}\frac{p_{\lambda}^{\prime}(|u^{(0)}|)}{|u^{(0)}|}(u^{2}-u^{(0)2}).
\]
Then given some initial estimate $\beta^{(0)}$, for $\beta\approx\beta^{(0)},$
we have
\begin{equation}
\sum_{j=1}^{M}p_{\lambda}\left(\frac{\|\beta_{[j]}\|_{2}}{\sqrt{T/M}}\right)\approx\frac{1}{2}\sum_{j=1}^{M}\frac{p_{\lambda}^{\prime}\left(\frac{\|\beta_{[j]}^{(0)}\|_{2}}{\sqrt{T/M}}\right)}{\frac{\|\beta_{[j]}^{(0)}\|_{2}}{\sqrt{T/M}}}\frac{\|\beta_{[j]}\|_{2}^{2}}{T/M}+G(\beta^{(0)}),\label{eq:LQA-Theoretical}
\end{equation}
where $G(\beta^{(0)})$ is a term defined as
\[
G(\beta^{(0)})\define\sum_{j=1}^{M}p_{\lambda}\left(\frac{\|\beta_{[j]}^{(0)}\|_{2}}{\sqrt{T/M}}\right)-\frac{1}{2}\sum_{j=1}^{M}p_{\lambda}^{\prime}\left(\frac{\|\beta_{[j]}^{(0)}\|_{2}}{\sqrt{T/M}}\right)\frac{\|\beta_{[j]}^{(0)}\|_{2}}{\sqrt{T/M}}.
\]
Let
\begin{equation}
\bdW^{(0)}=\frac{1}{2}\sum_{j=1}^{M}\left(\frac{p{}_{\lambda}^{\prime}(\|\beta_{[j]}^{(0)}\|_{2}\sqrt{M/T})}{\|\beta_{[j]}^{(0)}\|_{2}\sqrt{T/M}}\bdW_{j}\right).\label{eq:W0}
\end{equation}
Then we have
\begin{equation}
\frac{M}{T}\int_{0}^Tp_{\lambda}(|\beta(t)|)\diff t\approx\bdb^{T}\bdW^{(0)}\bdb+G(\beta^{(0)}).\label{eq:fscad-part}
\end{equation}

Now, we put (\ref{eq:ls-part}), (\ref{eq:roughness-part}) and (\ref{eq:fscad-part})
together to express $Q(\beta)$ in the matrix notation as
\begin{equation}
R(\bdb)\define\frac{1}{n}(\bdy-\bdU\bdb)^{T}(\bdy-\bdU\bdb)+\gamma\bdb^{T}\bdV\bdb+\bdb^{T}\bdW^{(0)}\bdb+G(\beta^{(0)}).\label{eq:objective-LQA-empirical}
\end{equation}
Thus, optimizing $Q(\beta)$ with respect to $\beta(t)$ is equivalent
to optimizing $R(\bdb)$ with respect to $\bdb$ in (\ref{eq:objective-LQA-empirical}).
Also note that the term $G(\beta^{(0)})$ does not depend on $\beta$
and hence has no impact on optimizing $R(\bdb)$. Differentiating
$R(\bdb)$ with respect to $\bdb$ and setting it to zero, we have
the following equation
\[
\bdU^{T}\bdU\bdb+n\gamma\bdV\bdb+n\bdW^{(0)}\bdb=\bdU^{T}\bdy
\]
with the solution
\[
\hat{\bdb}=\left(\bdU^{T}\bdU+n\gamma\bdV+n\bdW^{(0)}\right)^{-1}\bdU^{T}\bdy.
\]

As in \citet{Fan2001}, we repeat the above computation steps until
$\hat{\bdb}$ converges and in each iteration we manually shrink variables that are very small to zero. In summary, we have the following algorithm
to compute $\hat{\bdb}$ and obtain the estimator $\hat{\beta}(t)=\bdB^{\trans}(t)\hat{\bdb}$:
\begin{lyxlist}{00.00.0000}
\item [{$Step\:1:$}] Compute the initial estimate $\hat{\bdb}^{(0)}=\left(\bdU^{T}\bdU+n\gamma\bdV\right)^{-1}\bdU^{T}\bdy$;
\item [{$Step\:2:$}] Given $\hat{\bdb}^{(i)}$, compute $\bdW^{(i)}$
and $\hat{\bdb}^{(i+1)}=\left(\bdU^{T}\bdU+n\gamma\bdV+n\bdW^{(i)}\right)^{-1}\bdU^{T}\bdy$; if a variable is very small in magnitude and makes $\bdU^{T}\bdU+n\gamma\bdV+n\bdW^{(i)}$ almost singular or badly scaled so that inverting $\bdU^{T}\bdU+n\gamma\bdV+n\bdW^{(i)}$ is numerically unstable, then it is manually shrunk to zero;
\item [{$Step\:3:$}] Repeat Step 2 until the convergence of $\hat{\bdb}$
is reached.
\end{lyxlist}

To carry out the above algorithm, we need to evaluate the matrix $\mathbf{U}$ with elements $u_{ij}=\int_0^TX_i(t)B_j(t)\diff t$. When covariate functions $X_i(t)$ are not fully observed, approximation of elements $u_{ij}$ is mandatory. Provided that covariate functions $X_i(t)$ are observed at a regular and dense grid of time points $t_1,t_2,\ldots,t_K$ over $[0,T]$, we can compute $u_{ij}\approx\frac{1}{K}\sum_{k=1}^KX_i(t_k)B_j(t_k)$. An alternative approach is to fit each $X_i(t)$ using the spline regression method first. Then we compute $u_{ij}=\int_0^T\hat{X}_i(t)B_j(t)\diff t$, where $\hat{X}_i(t)$ is the fitted curve of $X_i(t)$ based on the observations $X_i(t_1),X_i(t_2),\ldots,X_i(t_K)$. When the grid is regular and dense, these two approaches yield  an almost identical approximate of $u_{ij}$. When $X_i(t)$ is observed on a sparse and irregular grid, we can first estimate $X_i(t)$ from the sparse measurements using some available methods such as functional principal component analysis (Yao et al. (2005)).

When $\mu\neq0$, we modify the matrices $\bdU$, $\bdV$ and $\bdW^{(i)}$
as follows to produce an estimate of $\hat{\mu}$ simultaneously:
1) add one more parameter $\mu$ to the topmost of the column vector
$\bdb$; 2) add one more column of all ones to the leftmost of $\bdU$; 3) 
add one column of all zeros to the leftmost of both $\bdV$ and $\bdW^{(i)}$,
and then add one more row of all zeros to the topmost of both $\bdV$
and $\bdW^{(i)}$. With these changes, the algorithm above can be
carried out to estimate $\beta(t)$ and $\mu$ simultaneously.

When using the smoothing spline method, the choice of $M$ is not crucial (\cite{Cardot2003}), as the roughness of the estimator is controlled by the roughness penalty, rather than the number of knots. This is demonstrated by a simulation study, which is detailed in Section A in the supplementary file. The simulation study also provides a practical guideline on choosing $M$: estimators with a larger $M$ perform better in identifying null region of $\beta(t)$, while those with a smaller $M$ perform better in prediction. Once $M$ is determined, we can perform a search over a grid of candidate values of $\gamma$ and $\lambda$ based on some popular selection criteria such as cross-validation, generalized cross-validation, BIC, AIC or RIC.

\subsection{Extension to Multiple Regressors\label{sub:Extension}}
Our method can also be extended to estimate the following functional linear model with multiple functional covariates:
\begin{equation}
Y_{i}=\mu+ \sum_{k=1}^K \int_0^T X_{ki}(t)\beta_k(t)\diff t+\varepsilon_{i},\label{eq:mflm}
\end{equation}
where each $Y_{i}$ is a scalar response, $\mu$ is the grand mean,
$X_{k1}(t), \ldots, X_{kn}(t)$ are independent realizations of an unknown
centered stochastic process $X_k(t)$ defined on the domain $[0,T]$,
and $\varepsilon_{i}\overset{i.i.d.}{\sim}N(0,\sigma_{\varepsilon}^{2})$, $i = 1, \ldots, n$. We estimate the coefficient functions, $\beta_k(t), k=1,\ldots,K$, and $\mu$ by minimizing
\begin{equation}
 Q(\boldsymbol{\beta}, \mu) = \frac{1}{n}\sum_{i=1}^{n}\left[Y_{i}-\mu-\sum_{k=1}^K\int_{0}^TX_{ki}(t)\beta_k(t)\diff t\right]^{2}+\sum_{k=1}^K\gamma_k\|\diffop^{m}\beta_k\|^{2}+\sum_{k=1}^K \frac{M}{T}\int_{0}^Tp_{\lambda_k}(|\beta_k(t)|)\diff t,\label{eq:mslos}
\end{equation}
where $\boldsymbol{\beta}(t)=(\beta_1(t),\beta_2(t),\ldots,\beta_K(t))^T$. The computational steps described in Section (\ref{sub:Locally-Sparse-Smooth})
are modified as follows to minimize $ Q(\boldsymbol{\beta}, \mu)$. Corresponding
to each regressor $X_{k\cdot}$, we compute the matrix $\bdU_{k}$, which has $n\times(M+d)$ entries $u_{ij}=\int_{0}^TX_{ki}(t)B_{j}(t)\diff t$
for $i=1,2,\ldots,n$, $j=1,2,\ldots,(M+d)$. Let
$\bdU=(\bdU_{1},\bdU_{2},\ldots,\bdU_{K})$ be the column catenation
of $\bdU_{1},\bdU_{2},\ldots,\bdU_{K}$, and correspondingly set $\bdV=\mathrm{diag}(\bdV_{1},\bdV_{2},\ldots,\bdV_{K})$, i.e., set $\bdV$ to be the matrix with blocks $\bdV_{1},\bdV_{2},\ldots,\bdV_{K}$ in its main diagonal and zeros elsewhere. For each $k$, we also have a matrix $\bdW^{(0,k)}$ which corresponds
to the matrix $\bdW^{(0)}$ in (\ref{eq:W0}). Then $\bdW^{(0)}$
in (\ref{eq:objective-LQA-empirical}) is replaced by $\bdW^{(0)}=\mathrm{diag}(\bdW^{(0,1)},\bdW^{(0,2)},\ldots,\bdW^{(0,K)})$.
After these modifications, the iterative algorithm described in Section
\ref{sub:Locally-Sparse-Smooth} can be carried out without any
change.

\section{Theoretical Properties\label{sec:asymptotic}}

Let us set up some notations first. We use $\Gamma_{n}$ to denote the empirical version of the covariance operator
$\Gamma$ of the random process $X$, and is defined by
\[
(\Gamma_{n}x)(v)\define\frac{1}{n}\sum_{i=1}^{n}\int_{D}X_{i}(v)X_{i}(u)x(u)\diff u.
\] 
For a function $f$, define a semi-norm $\|f\|_{\Gamma_n}=\langle \Gamma_nf,f\rangle^{1/2}$ and its population version $\|f\|_{\Gamma}=\langle \Gamma f,f\rangle^{1/2}$. Let $N(\beta)$ denote the null region of $\beta(t)$ and $S(\beta)$
denote the non-null region of $\beta(t)$, i.e., $N(\beta)=\{t\in [0,T]:\beta(t)=0\}$
and $S(\beta)=\{t\in [0,T]:\beta(t)\neq0\}$.

We also assume the following regularity conditions:
\begin{lyxlist}{00.00.0000}
\item [{(C1)}] $\|X\|_{2}$ is almost surely bounded, i.e., $\|X\|_{2}\leq c_{1}<\infty$
a.s for some constant $c_{1}>0$. Here, $\|X\|_{2}$ denotes $(\int_0^TX^2(t)\diff t)^{1/2}$.
\item [{(C2)}] $\beta(t)$ is in
the H\"{o}lder space $C^{p^{\prime},\nu}([0,T])$. That is, $|\beta^{(p^{\prime})}(u_{1})-\beta^{(p^{\prime})}(u_{2})|\leq c_{1}|u_{1}-u_{2}|^{\nu}$
for some constant $c_{1}>0$, integer $p^{\prime}$ and $\nu\in[0,1]$.
Let $p\define p^{\prime}+\nu$. Assume $3/2<p\leq d$, where $d$ is the
degree of the B-spline basis.
\item [{(C3)}] There exists a sequence of $\lambda_n$ such that $\lambda_n\converge 0$, such that $\sqrt{\int_{S(\beta)}p_{\lambda_n}^{\prime}\left(|\beta(t)|\right)^{2}\diff t}=O(n^{-1/4}M^{-1/2})$, and $\sqrt{\int_{S(\beta)}p_{\lambda_{n}}^{\prime\prime}(|\beta(t)|)^{2}\diff t}=o_{P}(M^{-3/2}n^{-1})$.
\end{lyxlist}
In the above, (C1) and (C2) are the same as (H1) and (H3) of \citet{Cardot2003}. Additionally, 
Assumption (C3) is analogous to regularity conditions ($\mathrm{B^{\prime}}$)
and ($\mathrm{C}^{\prime}$) in \citet{Fan2004} to ensure the unbiasedness
and guarantee that the penalty does not dominate the least squares. Intuitively, Assumption (C3) prevents the tail of $\beta(t)$ on its support from being too thin. If the tail is too thin, then the fSCAD penalty dominates the least squares loss.

Recall that $\gamma$ and $\lambda$ are tuning parameter varying
with $n$. To emphasize their dependency on the sample size $n$,
we denote them by $\gamma_{n}$ and $\lambda_{n}$, respectively.
In addition, we assume the following conditions
on choosing values of $M$, $\gamma_{n}$ and $\lambda_{n}$:
\begin{lyxlist}{00.00.0000}
\item [{(C4)}] $M/(n\gamma_n)=o(1)$ and $\gamma_n=o(n^{-1/2})$.
\end{lyxlist}

The theorem below establishes the existence of our estimator and its consistency.
\begin{theorem}\label{thm:local-solution}Under assumptions
(C1)-(C4), with probability tending to one, there exists a local minimizer
$(\hat{\beta},\hat{\mu})$ of (\ref{eq:slos}) such that $\|\hat{\beta}-\beta\|_{\Gamma_n}=o_P(1)$ and $\|\hat{\beta}-\beta\|_\Gamma=o_P(1)$.
\end{theorem}

We now introduce some notations for stating the next property. First,
recall that the support of $\beta(t)$ is defined as the closure
of $S(\beta)$. For any $\epsilon>0$ and a subset of the real line
$A$, the $\epsilon$-neighborhood of $A$, defined by $\{t\in [0,T]:\inf_{u\in A}|t-u|<\epsilon\}$,
is denoted by $A^{\epsilon}$. We also define $A^{-\epsilon}$ as $D-(D-A)^\epsilon$. Intuitively, $A^\epsilon$ extends $A$ a little bit, while $A^{-\epsilon}$ shrinks $A$ a little bit. The following theorem shows that, with probability tending to one, $\hat{\beta}(t)$ is collectively zero on the null
region except for those points on the boundaries.
\begin{theorem}\label{thm:oracle-property}Under conditions (C1)-(C4),
as $n\converge\infty$,
For every $t$ not in the support of $\beta(t)$,
we have $\hat{\beta}(t)=0$ with probability tending to one. Moreover,
for every $\epsilon>0$, we have $N^\circ(\beta)\subset N^{\epsilon}(\hat{\beta})$ and $S^{-\epsilon}(\beta)\subset S(\hat{\beta})$ 
with probability tending to one, where $N^\circ(\beta)$ denotes the interior of $N(\beta)$.
\end{theorem}

Proofs of the above theorem have been relegated to supplementary file.

\section{Simulation Studies\label{sec:Simulations}}

We conducted a simulation study  to investigate the numerical performance of our SLoS estimator. 
In this study, we compared our estimator to the following estimators:
1) the ordinary least squares (OLS) estimator, corresponding
to $\gamma=0$ and $\lambda=0$, 2) the smoothing spline (Smooth) estimator proposed in \citet{Cardot2003}, corresponding to $\lambda=0$, and 3) the oracle estimator. The
oracle estimator was computed by placing equally spaced knots only
on the non-null region $S(\beta)$ and was regularized by a roughness penalty. Notice that only the oracle estimator assumes the null region is known, while the other three estimators do not know the null region. It will be shown that our SLoS estimator performs much better than the OLS estimator and the smooth estimator, and the SLoS estimator performs almost as well as the oracle estimator, although our SLoS estimator does not know the null region in advance.  In our simulation study, we did not compare our method with \citet{Zhou2013} and \citet{James2009} because, the FLiRTI method proposed in \citet{James2009} is numerically unstable in our simulation settings, and the method of \citet{Zhou2013} is fairly difficult to implement and tune. This difficulty stems from the fact that, at each stage, several tuning parameters have to be chosen carefully.

In the comparison, we reported the integrated  squared error
($\mathrm{ISE}$) defined on the null region and the non-null region as follows
\begin{eqnarray}\label{Eqn:ISE}
\mathrm{ISE}_{0}=\frac{1}{\ell_{0}}\int_{N(\beta)}\left(\hat{\beta}(t)-\beta(t)\right)^{2}\diff t\qquad\quad \mathrm{ISE}_{1}=\frac{1}{\ell_{1}}\int_{S(\beta)}\left(\hat{\beta}(t)-\beta(t)\right)^{2}\diff t,
\end{eqnarray}
where $\ell_{0}$ was the length of the null region and $\ell_{1}$
was the length of the non-null region of $\beta(t)$. $\mathrm{ISE}_{0}$
and $\mathrm{ISE}_{1}$ measure the integrated squared error between an estimated coefficient function $\hat{\beta}(t)$
and the true function $\beta(t)$ on the null region and the non-null
region, respectively. We also assessed the  performance of estimators on prediction by the mean squared error (MSE) on predicting $y$ on a test dataset that is independent of the training dataset. The prediction mean sqaured error ($\mathrm{PMSE}$) is computed as follows:
\begin{eqnarray}\label{Eqn:PMSE}
\mathrm{PMSE}=\frac{1}{N}\sum_{(X,y)\in test}\left(y-\hat{\mu}-\int_0^TX(t)\hat{\beta}(t)\diff t\right)^2,
\end{eqnarray}
where $test$ denotes the test dataset, $N$ is the size of the test dataset, and $(\hat{\mu},\hat{\beta}(t))$ is the estimated intercept $\mu$ and coefficient function $\beta(t)$ from the training dataset. 

The SLoS estimator used cubic B-spline basis functions defined by $M+1$ equally spaced knots in [0,1]. For choosing $M$, we found that the empirical formula $M=\max\{50,[20n^{1/4}]\}$  works quite well in practice, where $[x]$ denotes the integer closest to $x$. According to this formula, $M\geq 50$. It made $M$ large enough to well approximate the fSCAD penalty. The values of $M$ in the following simulation studies were determined by the above formula. The tuning parameters $\gamma$ and $\lambda$
were chosen using BIC as the selection criterion by the procedure proposed in Section \ref{sub:Locally-Sparse-Smooth}. BIC was used because it encouraged sparser models. The smoothing spline method also used cubic B-spline basis functions, and used AIC to choose the smoothing parameter $\gamma$ and the number of knots that defined the B-spline basis, because AIC yielded better performance than cross-validation or BIC. The OLS estimator also used AIC to choose the order and the number of knots which defined the B-spline basis, because AIC also yielded better performance than cross-validation or BIC. The oracle estimator used cross-validation to choose the number of knots and the smoothing parameter $\gamma$, as cross-validation yielded better performance than AIC or BIC.

\begin{figure}[H]
\begin{minipage}[t]{0.49\columnwidth}%
\centering
\subfloat[]{
\includegraphics[scale=0.7]{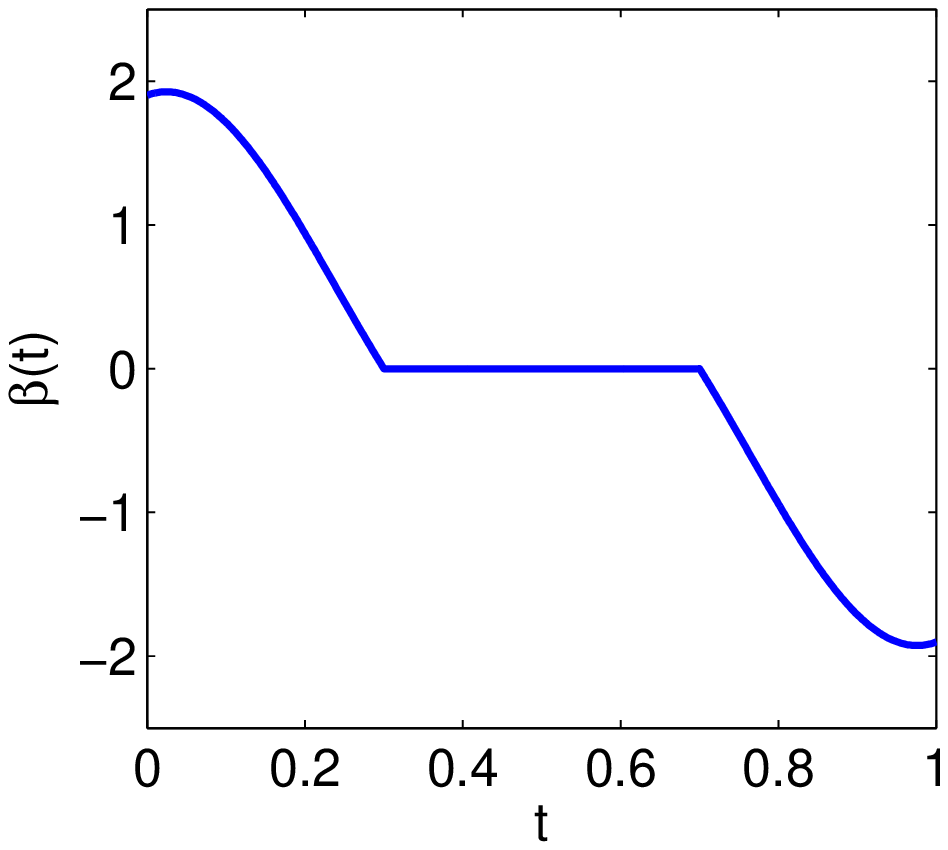}%
\label{subfig:case-II}}
\end{minipage}%
\begin{minipage}[t]{0.49\columnwidth}%
\centering
\subfloat[]{
\includegraphics[scale=0.7]{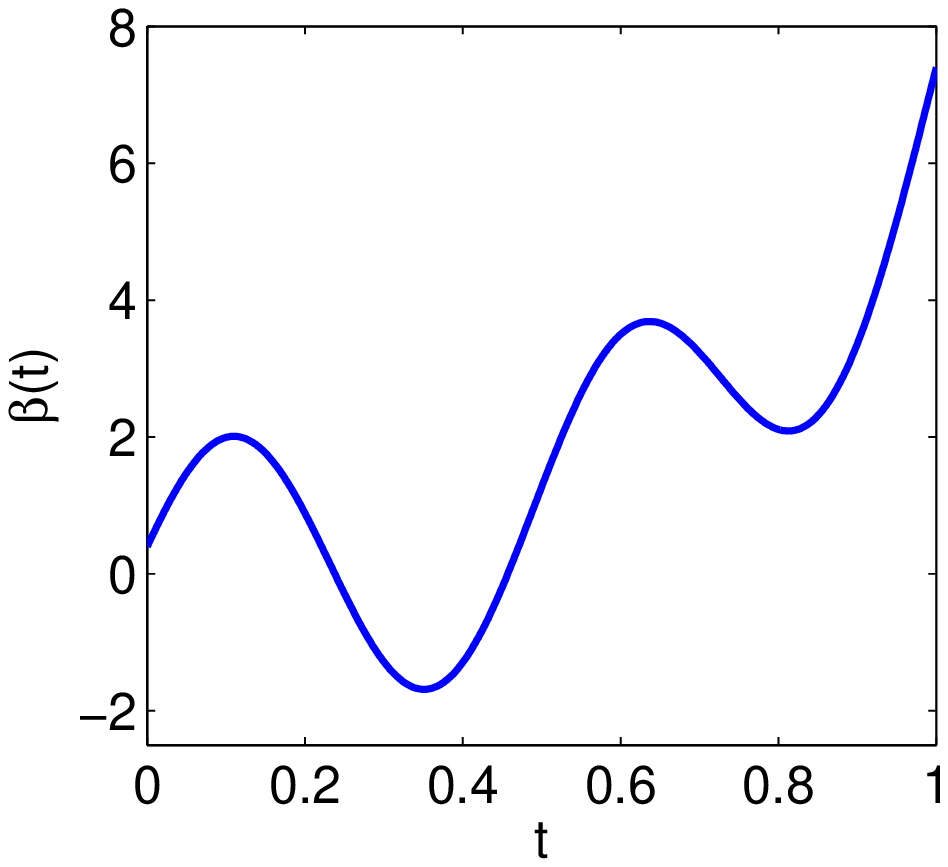}%
\label{subfig:case-III}}
\end{minipage}%
\caption{\protect\subref{subfig:case-II} the true $\beta(t)$ in Case I. \protect\subref{subfig:case-III} the true $\beta(t)$ in Case II.}
\label{fig:true-beta}
\end{figure}

We considered the functional linear model
\begin{eqnarray}\label{Eqn:FLMsimu}
Y_{i}=\mu + \int_{0}^{1}X_{i}(t)\beta(t)\diff t+\varepsilon_{i},
\end{eqnarray}
where $\varepsilon_i \sim N(0,\sigma_\varepsilon^2)$. Three different types of coefficient functions $\beta(t)$ were considered: 

\begin{lyxlist}{00.00.0000}
\item [{Case I:}] There is no signal, i.e., $\beta(t)\equiv0$.
\item [{Case II:}] There is a flat region of $\beta(t)$. We considered the following function  \[
\beta(t)=\begin{cases}
2(1-t)\sin(2\pi(t+0.2)) & 0\leq t\leq0.3,\\
0 & 0.3<t<0.7,\\
2t\sin(2\pi(t-0.2)) & 0.7\leq t\leq1,
\end{cases}
\]
which vanishes on $[0.3,0.7]$. The function is plotted in Figure \ref{fig:true-beta}\protect\subref{subfig:case-II}.
\item [{Case III:}] There is no flat region of $\beta(t)$. We considered $\beta(t)=4x^3+2\sin(4\pi t+0.2)$. It is plotted in Figure \ref{fig:true-beta}\protect\subref{subfig:case-III}. It was designed so that it has no flat region, but crosses zero twice. Therefore, this case does not favor the SLoS method.
\end{lyxlist}
The measurement error $\sigma_\varepsilon$ was set to 1 in Case I, and chosen so that the signal-to-noise ratio equals to 4 in the other two cases. The true $\mu$ was set to 1. 

The covariate functions $X_{i}(t)$ were generated based on the equation
$X_{i}(t)=\sum_{j=1}^{74}a_{ij}B_j(t)$, where the coefficients $a_{ij}$ were generated from the standard normal distribution, and $B_j(t)'s$ were B-spline basis functions defined by order 5 and 71 equally spaced  knots (the number 71 is randomly selected between 50 and 100). Using this model, for each sample size $n=150,450,1000$, we independently generated 100 datasets. For each dataset, we also generated a separate test dataset with the sample size 5000.

\begin{table}[htbp]
\begin{centering}
\begin{tabular}{lr@{\extracolsep{0pt}.}lr@{\extracolsep{0pt}.}lr@{\extracolsep{0pt}.}lr@{\extracolsep{0pt}.}l}
\hline 
 & \multicolumn{2}{c}{Oracle} & \multicolumn{2}{c}{OLS} & \multicolumn{2}{c}{Smooth} & \multicolumn{2}{c}{SLoS}\tabularnewline
\hline 
Case I $(\times10^{-2})$ & \multicolumn{2}{c}{} & \multicolumn{2}{c}{} & \multicolumn{2}{c}{} & \multicolumn{2}{c}{}\tabularnewline
~~~~~ $n=150$  & \multicolumn{2}{c}{-} & 2&10 (1.61) & 0&91 (0.92) & 0&23 (0.43)\tabularnewline
~~~~~ $n=450$ & \multicolumn{2}{c}{-} & 0&57 (0.46) & 0&29 (0.24) & 0&05 (0.11)\tabularnewline
~~~~~ $n=1000$ & \multicolumn{2}{c}{-} & 0&27 (0.22) & 0&12 (0.10) & 0&02 (0.03)\tabularnewline
Case II $(\times10^{-4})$ & \multicolumn{2}{c}{} & \multicolumn{2}{c}{} & \multicolumn{2}{c}{} & \multicolumn{2}{c}{}\tabularnewline
~~~~~ $n=150$  & 1&62 (0.89) & 3&89 (2.39) & 3&47 (2.18) & 2&11 (1.25)\tabularnewline
~~~~~ $n=450$ & 0&58 (0.32) & 1&28 (0.58) & 1&06 (0.40) & 0&72 (0.36)\tabularnewline
~~~~~ $n=1000$ & 0&28 (0.14) & 0&67 (0.31) & 0&57 (0.24) & 0&37 (0.16)\tabularnewline
Case III $(\times10^{-3})$ & \multicolumn{2}{c}{} & \multicolumn{2}{c}{} & \multicolumn{2}{c}{} & \multicolumn{2}{c}{}\tabularnewline
~~~~~ $n=150$ & \multicolumn{2}{c}{-} & 2&25 (1.42) & 2&12 (1.36) & 1&83 (0.78)\tabularnewline
~~~~~ $n=450$ & \multicolumn{2}{c}{-} & 0&64 (0.35) & 0&61 (0.28) & 0&56 (0.21)\tabularnewline
~~~~~ $n=1000$ & \multicolumn{2}{c}{-} & 0&35 (0.20) & 0&30 (0.15) & 0&28 (0.11)\tabularnewline
\hline 
\end{tabular}
\par\end{centering}
\caption{The prediction mean squared error ($\mathrm{PMSE}$) on test data using four methods: the oracle method, the ordinary least squares method (abbreviated as OLS), the smoothing spline method (abbreviated as Smooth), and our SLoS method. $\mathrm{PMSE}$ is defined in (\ref{Eqn:PMSE}). Each
entry is the Monte Carlo average of 100 simulation replicates. The corresponding Monte
Carlo standard deviation is included in parentheses.}
\label{tab:pmse}
\end{table}

\begin{table}[htbp]
\begin{centering}
\begin{tabular}{lcr@{\extracolsep{0pt}.}lr@{\extracolsep{0pt}.}l}
\hline 
 & OLS & \multicolumn{2}{c}{Smooth} & \multicolumn{2}{c}{SLoS}\tabularnewline
\hline 
Case I &  & \multicolumn{2}{c}{} & \multicolumn{2}{c}{}\tabularnewline
~~~~~ $n=150$ & 1.65 (1.67) & 0&57 (0.72) & 0&06 (0.31)\tabularnewline
~~~~~ $n=450$ & 0.42 (0.41) & 0&18 (0.18) & 0&01 (0.10)\tabularnewline
~~~~~ $n=1000$ & 0.20 (0.22) & 0&08 (0.08) & 0&00 (0.00)\tabularnewline
Case II $(\times10^{-3})$ &  & \multicolumn{2}{c}{} & \multicolumn{2}{c}{}\tabularnewline
~~~~~ $n=150$ & 19.75 (18.15) & 20&09 (18.20) & 0&15 (0.32)\tabularnewline
~~~~~ $n=450$ & 7.34 (4.62) & 6&37 (3.31) & 0&04 (0.09)\tabularnewline
~~~~~ $n=1000$ & 4.23 (3.10) & 4&00 (2.58) & 0&01 (0.03)\tabularnewline
\hline 
\end{tabular}

\par\end{centering}
\caption{The integrated squared error, $\mathrm{ISE}_0$, defined on the null region for the estimator using three methods: the ordinary least squares method (abbreviated as OLS), the smoothing spline method (abbreviated as Smooth), and our SLoS method. $\mathrm{ISE}_0$ is defined in (\ref{Eqn:ISE}). Each entry is the Monte Carlo average of 100 simulation replicates. The corresponding Monte
Carlo standard deviation is included in parentheses.}
\label{tab:ise0}
\end{table}

\begin{table}[htbp]
\begin{centering}
\begin{tabular}{lr@{\extracolsep{0pt}.}lr@{\extracolsep{0pt}.}lr@{\extracolsep{0pt}.}lr@{\extracolsep{0pt}.}l}
\hline 
 & \multicolumn{2}{c}{Oracle} & \multicolumn{2}{c}{OLS} & \multicolumn{2}{c}{Smooth} & \multicolumn{2}{c}{SLoS}\tabularnewline
\hline 
Case II $(\times10^{-2})$ & \multicolumn{2}{c}{} & \multicolumn{2}{c}{} & \multicolumn{2}{c}{} & \multicolumn{2}{c}{}\tabularnewline
~~~~~ $n=150$  & 1&93 (1.18) & 4&13 (3.56) & 3&11 (2.68) & 2&51 (1.60)\tabularnewline
~~~~~ $n=450$  & 0&73 (0.44) & 1&23 (0.88) & 0&88 (0.47) & 0&86 (0.44)\tabularnewline
~~~~~ $n=1000$  & 0&37 (0.18) & 0&66 (0.39) & 0&46 (0.22) & 0&46 (0.19)\tabularnewline
Case III $(\times10^{-1})$ & \multicolumn{2}{c}{} & \multicolumn{2}{c}{} & \multicolumn{2}{c}{} & \multicolumn{2}{c}{}\tabularnewline
~~~~~ $n=150$ & \multicolumn{2}{c}{-} & 1&88 (1.55) & 1&64 (1.35) & 1&32 (0.61)\tabularnewline
~~~~~ $n=450$ & \multicolumn{2}{c}{-} & 0&52 (0.36) & 0&46 (0.25) & 0&41 (0.17)\tabularnewline
~~~~~ $n=1000$ & \multicolumn{2}{c}{-} & 0&28 (0.20) & 0&22 (0.12) & 0&20 (0.08)\tabularnewline
\hline 
\end{tabular}
\par\end{centering}
\caption{The integrated squared error, $\mathrm{ISE}_1$, defined on the non-null region for the estimator using four methods: the oracle method, the ordinary least squares method (abbreviated as OLS), the smoothing spline method (abbreviated as Smooth), and our SLoS method. $\mathrm{ISE}_1$ is defined in (\ref{Eqn:ISE}).  
Each
entry is the Monte Carlo average of 100 simulation replicates. The corresponding Monte
Carlo standard deviation is included in parentheses.}
\label{tab:ise1}
\end{table}

The simulation results are shown in Table \ref{tab:pmse}, \ref{tab:ise0} and \ref{tab:ise1}. They clearly show that the SLoS method wins under each evaluation criterion in all cases except the
oracle procedure. From Table \ref{tab:pmse}, we can see that, in terms of prediction, our method performs significantly better than the OLS and Smooth estimators in Case I and II. Surprisingly, it also performs slightly better than the other two methods in Case III, although this case does not favor the SLoS estimator. We investigated this phenomenon, and found that the shrinkage effect of fSCAD penalty reduced the variability of SLoS estimator without causing significant bias. The reduction of variability is more significant when sample size is small. Therefore, our estimator yielded a better variance-bias tradeoff on the performance of prediction. We also noticed that the oracle procedure only performed slightly better than the SLoS method. 

In terms of estimation of $\beta(t)$, Table \ref{tab:ise0} shows that the SLoS method yielded a much lower integrated squared error on the null region, and Table \ref{tab:ise1} shows that it has a similar performance to the Smooth estimator on the non-null region. Even in Case III where the setting does not favor the SLoS method, the SLoS estimator still slightly outperforms the OLS and Smooth estimators. In summary, the SLoS method stands out in both predicting response variable and estimating the coefficient function $\beta(t)$.

\begin{figure}
\begin{minipage}[t]{0.47\columnwidth}%
\centering
\subfloat[]{
\includegraphics[scale=0.7]{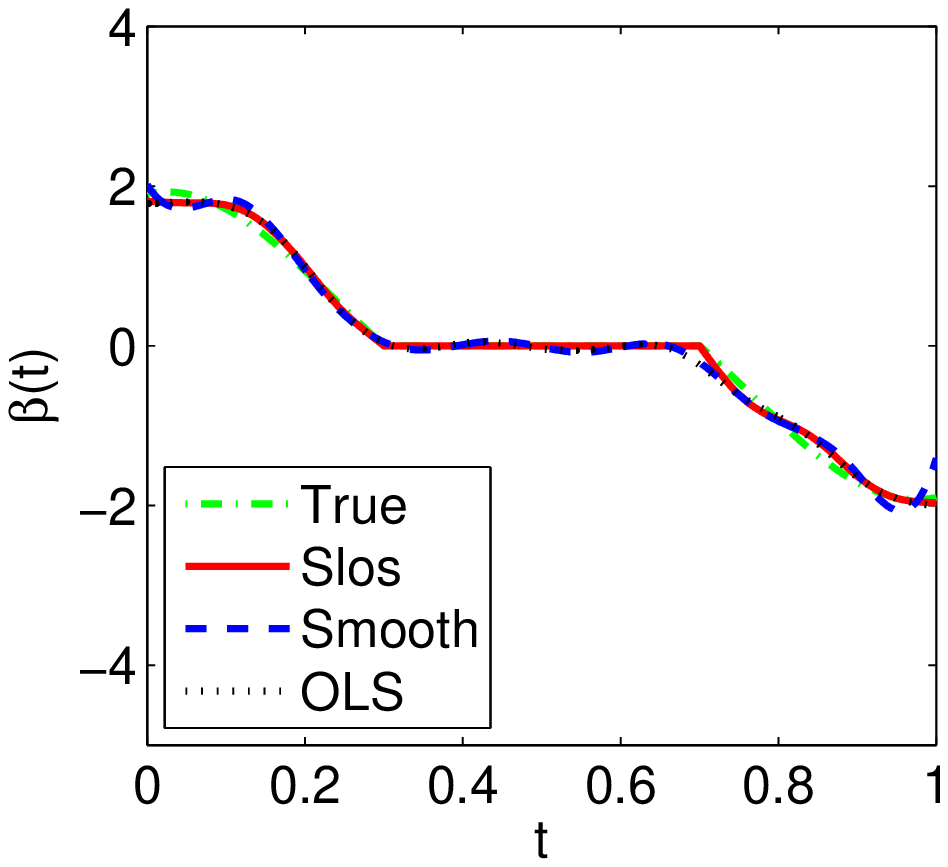}
\label{subfig:rand-s1}}%
\end{minipage}%
\begin{minipage}[t]{0.49\columnwidth}%
\centering
\subfloat[]{
\includegraphics[scale=0.7]{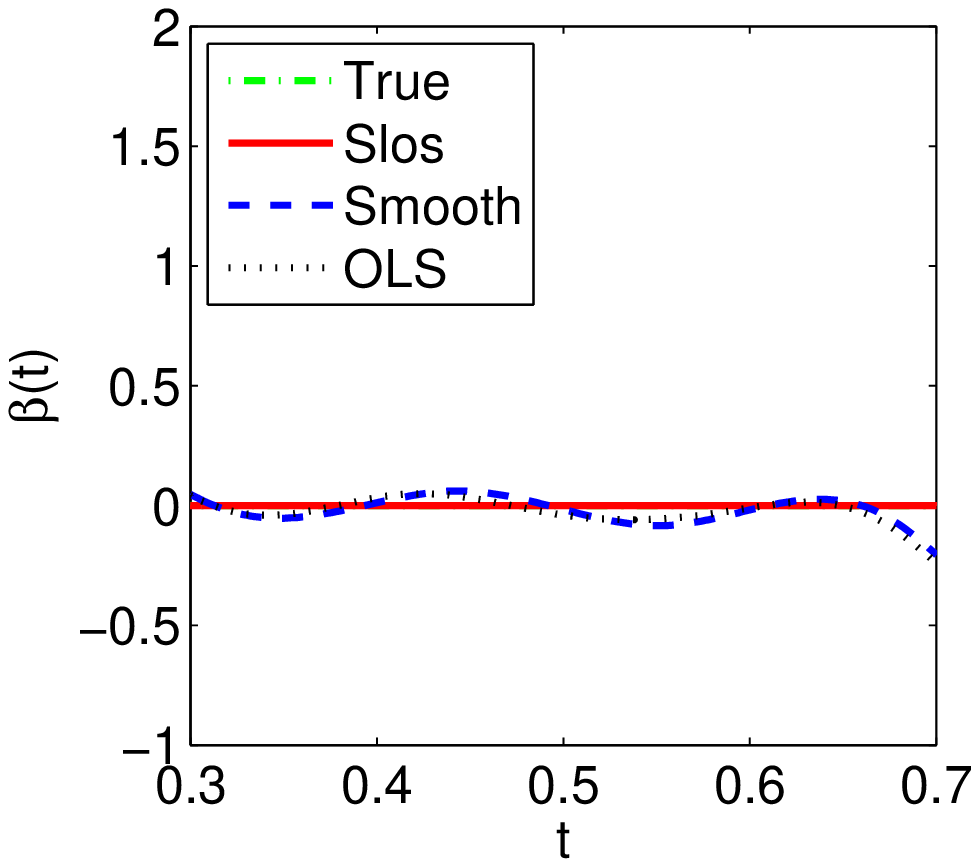}%
\label{subfig:rand-s1-n1}}%
\end{minipage}%
\caption[]{\protect\subref{subfig:rand-s1} shows the estimated $\hat{\beta}(t)$ in a simulation replicate chosen randomly in
the simulation. \protect\subref{subfig:rand-s1-n1} is the same $\hat{\beta}(t)$ in \protect\subref{subfig:rand-s1}, while it is zoomed in on the null region $[0.3,0.7]$ of $\beta(t)$. }
\label{Fig:one-random-run-study1}
\end{figure}

We also evaluate the ability of SLoS method to identify the null subregions. Figure \ref{Fig:one-random-run-study1}\protect\subref{subfig:rand-s1} displays the estimated $\hat{\beta}(t)$ by SLoS, OLS and Smooth method in one random simulation replicate, as well as the true $\beta(t)$. The estimated $\hat{\beta}(t)$ by the oracle procedure is not displayed in Figure \ref{Fig:one-random-run-study1}\protect\subref{subfig:rand-s1}, as it is almost identical to the true $\beta(t)$ and is hard to be distinguished in the graph. Figure \ref{Fig:one-random-run-study1}\protect\subref{subfig:rand-s1} shows that all three methods produce a very good estimate of $\beta(t)$. To get a close look at the null sug-region, Figure \ref{Fig:one-random-run-study1}\protect\subref{subfig:rand-s1-n1} displays the same functions of Figure \ref{Fig:one-random-run-study1}\protect\subref{subfig:rand-s1} only in the null region $[0.3,0.7]$. In this figure, we can see that, for the SLoS estimator, it is identically zero. However, the estimated $\hat{\beta}(t)$ produced by the other two methods are non-zero almost everywhere in the null region $[0.3,0.7]$. We also observe similar results in most of other simulation replicates.

Finally, we quantify the ability to identify the null subregions using the average of proportions of null subregions that are correctly identified by the SLoS estimator, which is computed as follows. First, for each run, we 
compute the value of $\hat{\beta}(t)$ on a sequence of dense and equally-spaced points in the null subregion. For example, in Case I, the sequence is taken to be $0,0.001,0.002,\ldots,0.999,1$, and in Case II, it is taken to be $0.3,0.301,0.302,\ldots,0.7$. Then we compute the proportion of the points at which $\hat{\beta}(t)$ is zero among all points in the sequence. Finally, we average the calculated proportions from 100 runs. The result is summarized in Table \ref{tab:null}. It shows that the performance of our estimator on identifying null region is quite remarkable: in Case I (Case II, respectively), on average more than 95\% (92\%) of the null region was correctly identified when $n=150$, and this number increases to 100\% (95\%) when $n=1000$. This also numerically verifies the conclusion in Theorem \ref{thm:oracle-property}.

\begin{table}[htbp]
\begin{centering}

\begin{tabular}{cccc}
\hline 
 & $n=150$ & $n=450$ & $n=1000$\tabularnewline
\hline 
Case I (\%) & 98.40 (6.14) & 99.70 (1.75) & 100 (0.00)\tabularnewline
Case II (\%) & 92.20 (5.47) & 93.41 (3.27) & 95.01 (2.04)\tabularnewline
\hline 
\end{tabular}
\par\end{centering}
\caption{The proportion of null region that is correctly identified by SLoS estimator on test data. The proportion of null region identified by SLoS estimator is computed as the proportion
of points where $\hat{\beta}(t)$ is zero on a sequence of dense and equally spaced points on null region. 
Each
entry is the Monte Carlo average of 100 simulation replicates. The corresponding Monte
Carlo standard deviation is included in parentheses.}
\label{tab:null}
\end{table}

\section{Applications}\label{sec:Applications}
\subsection{Canadian Weather Data}

\begin{figure}[htbp]
\begin{centering}
\includegraphics[scale=0.8]{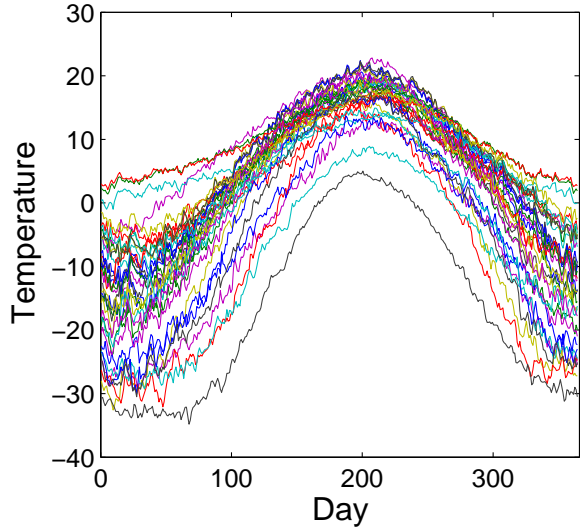}
\par\end{centering}

\caption{Daily mean temperature recorded at 35 Canadian cities in one year.}
\label{fig:weather-temp-curves}

\end{figure}

\cite{Ramsay2005} studied the relationship between the daily mean temperature
curves and the annual rainfall by analyzing the weather data of 35 Canadian cities. This relationship is modelled using the functional linear regression model (\ref{eq:flr-model}), where  the daily mean temperature curve is the functional covariate, and the logarithm of the annual rainfall is
the response variable.  Figure \ref{fig:weather-temp-curves} displays the daily mean temperature recorded at 35 Canadian cities
in one year. 

We applied our SLoS approach to estimate the functional linear regression model (\ref{eq:flr-model}) from the Canadian weather
data. The coefficient
function $\beta(t)$ was estimated as a linear combination of cubic B-spline basis functions defined
on $[0,365]$ with $M=50$. Here, the value of $M$ was also set by the empirical formula $M=\max\{50,[20n^{1/4}]\}$ that is introduced in the previous section.  To
respect the periodic nature of temperature cycles,  we restricted that $\beta(0) = \beta(365)$. The tuning parameters
were chosen by the the procedure proposed in
Section \ref{sub:Locally-Sparse-Smooth}.

\begin{figure}
\begin{centering}
\includegraphics[scale=0.7]{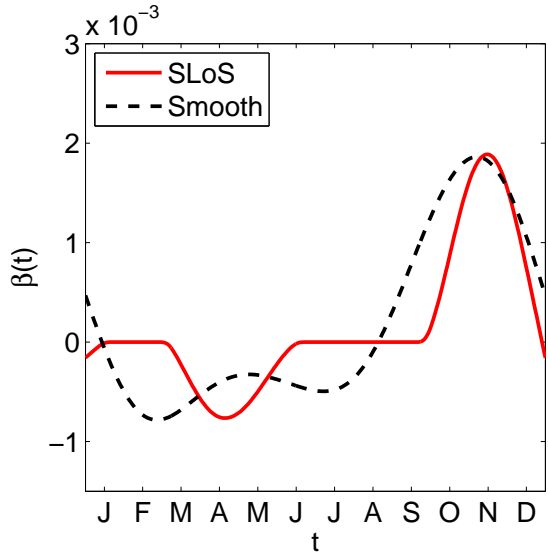}
\par\end{centering}

\caption{The estimated coeffient function $\hat{\beta}(t)$ for the functional linear model  (\ref{eq:flr-model}) from the Canadian weather data using our SLoS method (solid line) and the smoothing spline method (dotted line).}
\label{fig:weather-result}
\end{figure}

Figure \ref{fig:weather-result}
shows the temperature coefficient function estimated by our SLoS method
and the smoothing spline method \citep{Cardot2003}.
It shows that the estimated coefficient function $\hat{\beta}(t)$ is zero
roughly in January and summer months (late June, July, August and
 early September). This suggests the temperature in summer has no significant effect
on the annual rainfall. Our method also produced a smooth
estimate of $\beta(t)$ on the non-null subregions. It indicates
that the temperature in fall months (late October, November and early December) has a
significant contribution to the annual rainfall. These results are consistent with the results discovered in the previous research
on this data \citep{James2009}.

\begin{figure}[H]
\begin{centering}

\includegraphics[scale=0.65]{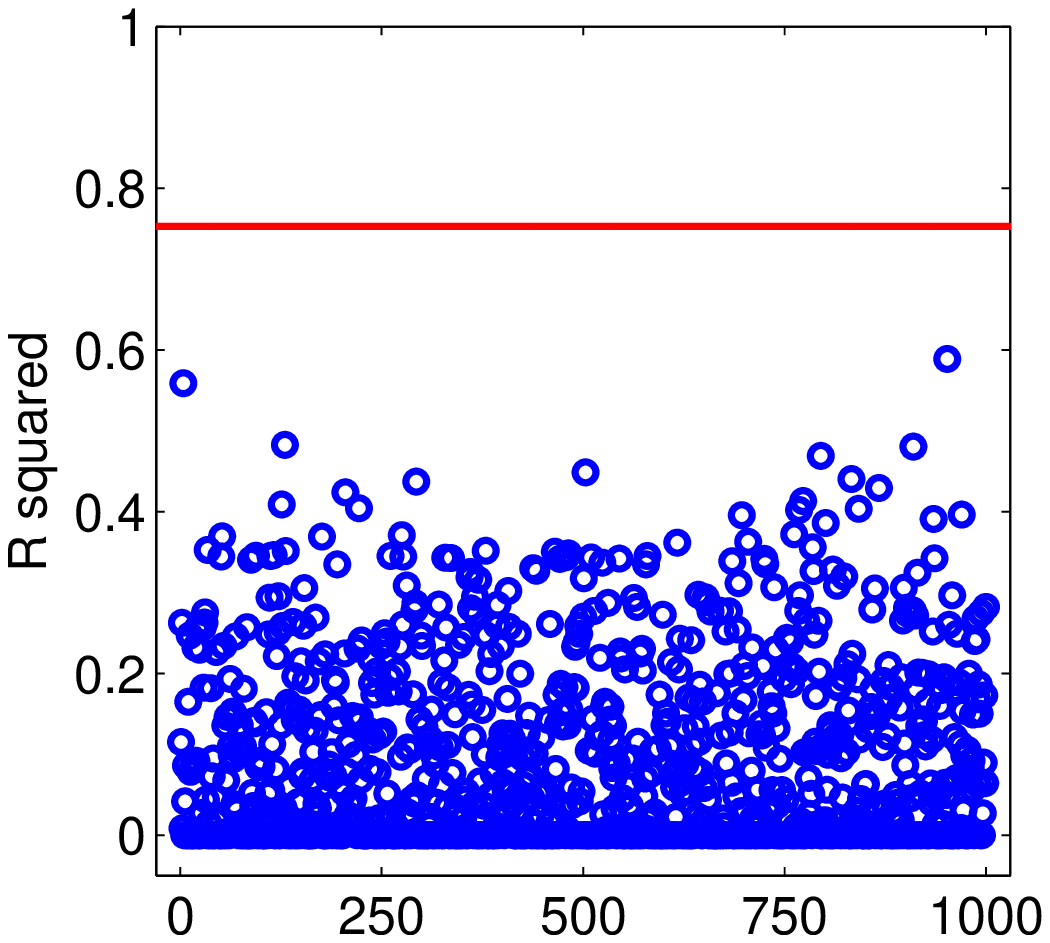}%

\par\end{centering}
\caption{$R^2$  from permuting the response variable 1000 times. The solid line represents the observed $R^2$ of SLoS estimate from the true data.}

\label{fig:weather-conf-test}
\end{figure}

Figure \ref{fig:weather-conf-test} shows the results of a permutation test for SLoS estimator on this data. The solid line represents the $R^2$ for the SLoS estimator, which is 0.73. The responses are random permuted 1000 times and the new $R^2$ are plotted in Figure \ref{fig:weather-conf-test}, indicated by circles. All 1000 permuted $R^2$ are under 0.73. This provides strong evidence of the relationship between temperature and annual rainfall discovered by our SLoS method.

\subsection{Spectrometric Data}

\begin{figure}[htbp]
\begin{centering}
\includegraphics[scale=0.8]{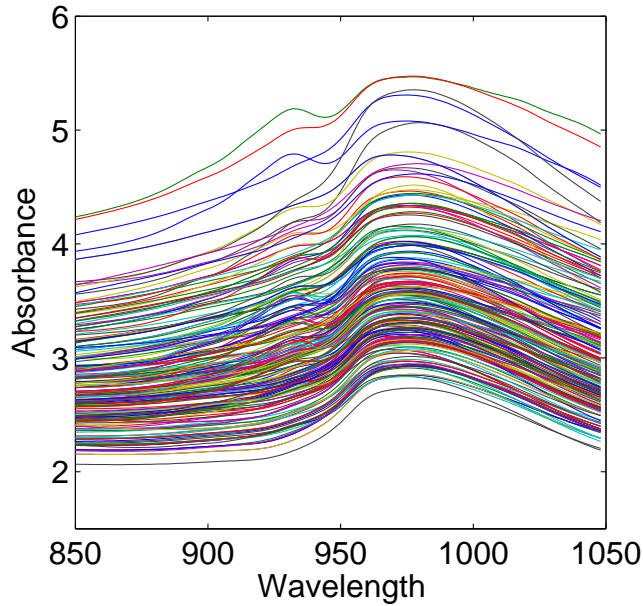}
\par\end{centering}

\caption{Spectrometric curves measured at 100 wavelengths among 850 - 1050 nm.}
\label{fig:spectrometric-curves}

\end{figure}

The data is available at http://www.math.univ-toulouse.fr/ferraty/SOFTWARES/NPFDA/npfda-datasets.html, which originates from Tecator dataset (http://lib.stat.cmu.edu/datasets/tecator). There are in total 215 samples. Each sample contains finely chopped pure meat, and a spectrometric curve of spectra of absorbances measured at 100 wavelengths among 850 - 1050 nm is recorded. Figure \ref{fig:spectrometric-curves} displays the 215 curves in the dataset.  At the same time, the fat content, measured in percent, is also determined by analytic chemistry. The task is to predict the fat content based on the spectrometric curve. Our interest is to investigate what range of spectra that has predicting power on fat content and what range that does not. 

We applied our SLoS approach to estimate the functional linear regression model (\ref{eq:flr-model}) from the spectrometric
data. The coefficient
function $\beta(t)$ was estimated as a linear combination of cubic B-spline basis functions defined
on $[850,1050]$ with $M=77$. Here, the value of $M$ was also set by the empirical formula $M=\max\{50,[20n^{1/4}]\}$ that is introduced in the previous section. The tuning parameters
were chosen by the the procedure proposed in
Section \ref{sub:Locally-Sparse-Smooth}.

\begin{figure}
\begin{centering}
\includegraphics[scale=0.7]{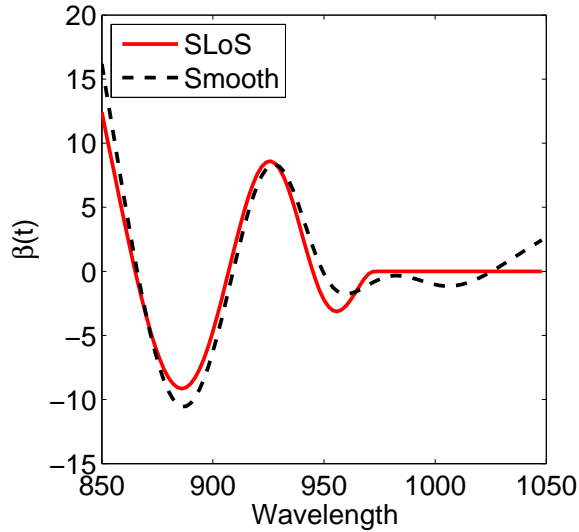}
\par\end{centering}

\caption{The estimated coeffient function $\hat{\beta}(t)$ for the functional linear model  (\ref{eq:flr-model}) from the spectrometric data using our SLoS method (solid line) and the smoothing spline method (dotted line).}
\label{fig:spectrometric-result}
\end{figure}

\begin{figure}[H]
\begin{centering}

\includegraphics[scale=0.65]{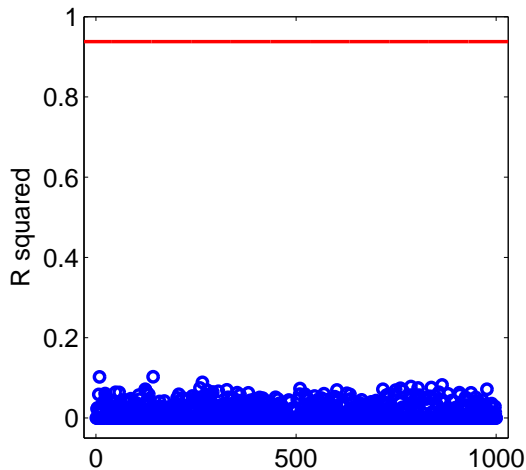}%
\par\end{centering}

\caption{$R^2$  from permuting the response variable 1000 times in spectrometric data. The solid line represents the observed $R^2$ of SLoS estimate from the true data.}

\label{fig:spectrometric-conf-test}
\end{figure}

Figure \ref{fig:spectrometric-result}
shows the coefficient function estimated by our SLoS method
and the smoothing spline method \citep{Cardot2003}.
It shows that the estimated coefficient function $\hat{\beta}(t)$ is zero
roughly on [970,1050]. This suggests the high end of spectrum has no predicting power on fat content. Our method also produced a smooth
estimate of $\beta(t)$ on the non-null subregions. It indicates
spectrum channels with wavelength among 850 - 970 nm have a significant contribution to the prediction power of the fat content.  Figure \ref{fig:spectrometric-conf-test} shows the results of a permutation test for SLoS estimator on this data. The solid line represents the $R^2$ for the SLoS estimator, which is 0.94. The responses are randomly permuted 1000 times and the new $R^2$ are plotted in Figure \ref{fig:spectrometric-conf-test}, indicated by circles. All 1000 permuted $R^2$ are under 0.2. This provides strong evidence of the relationship between spectrum and fat content discovered by our SLoS method.

\section{Concluding Remarks\label{sec:Conclusions}}

When a large number of predictors are involved in a regression problem,
parsimonious models via SCAD or other shrinkage regularization methods
have been proven to have less variability and better interpretability. In this paper, 
we propose the functional SCAD regularization method which extends
the ordinary SCAD to the functional setting. The functional SCAD,
when combined with the penalized B-spline expansion method, yields
a smooth and locally sparse (SLoS) estimator of the coefficient function
in functional linear regression.

The SLoS procedure is a combination
of three techniques: 1) the functional SCAD that is responsible for
identifying the null subregions of the coefficient function while
at the same time avoiding over shrinking the non-zero values, 2) the
B-spline basis expansion that is used to practically compute the SLoS
estimator efficiently thanks to its compact support property, and 3) 
the roughness regularization that assures the smoothness of our estimator
even when a large number of knots are used to define the B-spline
basis. Therefore, our method is able to accurately identify the null
region and simultaneously produce a smooth estimator on the non-null
region. Comparing to existing
methods in the literature, our estimation procedure is more theoretically sound and is computationally simpler. The simulation studies show
that our estimator has superior numerical performance and finite-sample
properties. Furthermore, the applications on two real datasets demonstrate
the practical merit of our method.

In summary, while our work focused on the functional linear regression, it is important to recognize that the functional SCAD is in fact a very general regularization technique which can  be applied in many other domains in functional data analysis. For example, it may be used in a spline
smoothing problem to obtain a smooth and locally sparse estimator
of an unknown curve.  It can also be used in functional principle component analysis to produce smooth and locally sparse functional principal components. Both of these problems are suggested starting points for future research.

\bibliographystyle{apalike}
\bibliography{LocallySparseFLR,VariableSelection,CLT}

\end{document}